\documentclass{article}

\usepackage{delarray,verbatim,enumerate,a4wide}
\usepackage{amsmath,amsthm,amstext,amsbsy,amssymb,amsfonts,amscd}
\usepackage[all]{xy}
\usepackage[frenchb]{babel}
\usepackage[T1]{fontenc}
\usepackage[utf8]{inputenc}

\usepackage{scalerel}

\usepackage[small]{titlesec}

\usepackage{color}
\definecolor{vert}{rgb}{0.1,0.4,0.2}
\usepackage[colorlinks=true,linkcolor=blue,citecolor=vert,linktocpage]{hyperref}

\usepackage[mathscr]{euscript}

\newtheorem{Th}{Théorème}[]
\newtheorem{Lem}[Th]{Lemme}
\newtheorem{Prop}[Th]{Proposition}
\newtheorem{Cor}[Th]{Corollaire}

\newtheorem{Sco}[Th]{Scolie}

\def\Preuve{\noindent {\it Preuve.~}}

\def\Remarque{\smallskip\noindent {\it Remarque.~}}

\def\Exemple{\smallskip\noindent {\it Exemple.~}}

\font\teneufm=eufm10
\font\seveneufm=eufm7
\font\fiveeufm=eufm5
\newfam\gothfam
\textfont\gothfam=\teneufm \scriptfont\gothfam=\seveneufm
\scriptscriptfont\gothfam=\fiveeufm

		\def\QQ{\mathbb Q}	
\def\NN{\mathbb N}	\def\ZZ{\mathbb Z}		
\def\F2{\mathbb{F}_2}	\def\Z2{\mathbb{Z}_2}		
\def\Zl{{\mathbb{Z}_\ell}} 	\def\Ql{\mathbb{Q}_\ell}		

 				\def\U{\mathcal  U}		\def\F{\mathcal  F}
  		\def\C{\mathcal  C}		\def\R{\mathcal  R}		
 	  	\def\Cl{\mathcal  C\!\ell}	
\def\E{\mathcal  E}		\def\T{\mathcal  T}			
\def\N{\mathcal N}		\def\Z{\mathcal Z}		\def\Cap{{\mathcal C\!ap}}

	\def\p{{\mathfrak p}}			
		\def\l{{\mathfrak l}}

\def\wi{\widetilde}

\def\Gal{\operatorname{Gal}}		\def\Rad{\operatorname{Rad}}
\def\Ker{\operatorname{Ker}}		

\def\ph{{\phantom{*}}}

\newcommand\scale[2]{\vstretch{#1}{\hstretch{#1}{#2}}}

\makeatletter
\newcommand*\wt[2][0.2ex]{%
        \begingroup
        \mathchoice{\wt@helper{#1}{#2}{\displaystyle}{\textfont}}
                   {\wt@helper{#1}{#2}{\textstyle}{\textfont}}
                   {\wt@helper{#1}{#2}{\scriptstyle}{\scriptfont}}
                   {\wt@helper{#1}{#2}{\scriptscriptstyle}{\scriptscriptfont}}%
        \endgroup
        #2%
}
\newcommand*\wt@helper[4]{%
        \def\currentfont{\the#41}%
        \def\currentskewchar{\char\the\skewchar\currentfont}%
        \setbox\tw@\hbox{\currentfont$#2$\currentskewchar}%
        \dimen@ii\wd\tw@
        \setbox\tw@\hbox{\currentfont$#2${}\currentskewchar}%
        \advance\dimen@ii-\wd\tw@
        \rlap{\raisebox{-#1}{$\m@th#3\kern\dimen@ii\widetilde{\phantom{#2}}$}}%
}
\makeatother

\def\wE{\,\wt[0.1ex]{\!\mathcal E}}	\def\we{\wt[0.3ex]{\mathfrak E}}	\def\wU{\wt[0.2ex]{\mathcal U}}
	\def\wCl{\wt[0.1ex]{\mathcal C\!\ell}} 
				
\def\wCap{\wt[0.2ex]{\Cap}}
\newcommand\si[1]{\scale{.8}{#1}}
\def\%{{\scale{.8}{\infty}}}

\newcolumntype{x}[1]{>{\centering\hspace{0pt}}p{#1}}

\begin{document}

\title{\Large\bf Normes universelles et conjecture de Greenberg}

\author{ Jean-François {\sc Jaulent} }
\date{}
\maketitle
\bigskip

{\small
\noindent{\bf Résumé.} Nous étudions le groupe des normes universelles attaché à la $\Zl$-tour cyclotomique d'un corps de nombres totalement réel à la lumière de la conjecture de Geenberg sur la trivialité des invariants d'Iwasawa.}

\

{\small
\noindent{\bf Abstract.} We investigate the group of universal norms attached to the cyclotomic $\Zl$-tower of a totally real number field in connection with  Grenberg's conjecture on Iwasawa invariants of such a field.}
\bigskip

\tableofcontents
\section*{Introduction}
\addcontentsline{toc}{section}{Introduction}

Étant donnés un corps de nombres $K$ et un nombre premier $\ell$ arbitraires, le pro-$\ell$-groupe des normes universelles de $K$, étudié par Greither dans \cite{Grt} et dont il est question ici, est le groupe universel pour le foncteur norme attaché à la $\Zl$-extension cyclotomique $K_{\si{\infty}}=\bigcup_{n\in\NN} K_n$ de $K$, relatif aux $\ell$-adifiés $\,\E'_{K_n}=\Zl\otimes_{\ZZ}E'_{K_n}$ des groupes de $\ell$-unités $E'_{K_n}$ des divers étages $K_n$.\smallskip

C'est tout simplement l'intersection des groupes de normes $\,\N_K=\bigcap_{n\in\NN}N_{K_n/K}(\E'_{K_n})$ et sa définition est irréductiblement globale.\smallskip

Il se présente ainsi comme un sous-groupe naturel du pro-$\ell$-groupe des normes cyclotomiques $\,\wE_K=\bigcap_{n\in\NN}N_{K_n/K}(\R_{K_n})$, défini, lui, comme intersection des groupes de normes attachés aux $\ell$-adifiés $\R_{K_n}=\Zl\otimes_{\ZZ}K_n^\times$ des groupes multiplicatifs des $K_n$, et qui peut donc être décrit localement en vertu du principe de Hasse (cf. e.g. \cite{J55}, App.). De fait, $\,\wE_K$ n'est rien d'autre que le pro-$\ell$-groupe des unités logarithmiques introduit dans \cite{J28}.\smallskip

Le but de la présente note est de comparer ces deux groupes qui jouent un rôle central dans plusieurs questions de la Théorie d'Iwasawa et, plus précisément, d'étudier le quotient $\,\wE_K/\N_K$ à la lumière de la conjecture de Greenberg (cf. \cite{Grb}).\smallskip

Notre point de départ est le Théorème \ref{TP} infra, qui relie l'indice $(\,\wE_K^\ph : \N_K\,)$ au cardinal du $\ell$-groupe des classes logarithmiques $\,\wCl_K$ du corps considéré.

\newpage
\section{Énoncé de la conjecture en termes de normes universelles}

La conjecture de Greenberg pour un nombre premier $\ell$ postule que les $\ell$-groupes de classes d'idéaux $\,\Cl_{K_n}$ attachés aux étages finis $K_n$ de la $\Zl$-extension cyclotomique d'un corps de nombres totalement réel $K$ sont d'ordre borné indépendamment de $n$.\smallskip

Sous la conjecture de Leopoldt, donc inconditionnellement pour $K$ abélien sur $\QQ$, il est bien connu que les sous-groupes sauvages $\,\Cl_{K_n}^{\,\si{[\ell]}}$, i.e. les sous-groupes respectifs des $\Cl_{K_n}^\ph$ construits sur les places au-dessus de $\ell$, sont d'ordre borné; de sorte que la conjecture de Greenberg revient à postuler qu'il en est de même des quotients $\,\Cl'_{K_n}=\Cl_{K_n}^\ph/\Cl_{K_n}^{\,\si{[\ell]}}$, i.e. des $\ell$-groupes de $\ell$-classes.\smallskip

En d'autres termes la Conjecture affirme que la limite projective pour les applications normes\footnote{Ce $\Lambda$-module, dit de Kuz'min-Tate, (ou de Tate par Kuz'min \cite{Kuz})  est noté $\,\T_{K_{\si{\infty}}}$ dans \cite{J55} et $\,\C'_{K_{\si{\infty}}}$ dans \cite{J60}.}\smallskip

\centerline{$\T^\ph_{K_{\si{\infty}}} = \,\C'_{K_{\si{\infty}}} = \varprojlim \,\Cl'_{K_n}$}\smallskip

\noindent regardée comme module sur  l'algèbre d'Iwasawa $\Lambda=\Zl[[\gamma-1]]$ construite sur un générateur topologique $\gamma$ du groupe procyclique $\Gamma=\Gal(K_{\si{\infty}}/K)$ est pseudo-nulle.\smallskip

Par un argument classique de Théorie d'Iwasawa (cf. e.g. \cite{J60}, Lem. 1), cela revient à postuler que le sous-module des points fixes $\,\T{\!}_{K_{\si{\infty}}}^{\;\Gamma}$ et le quotient de co-points fixes ${}^\Gamma\T_{K_{\si{\infty}}}$ sont deux $\Zl$-modules finis et de même ordre.
Or, ces deux modules ont une interprétation simple:\smallskip
\begin{itemize}
\item Le quotient des genres ${}^\Gamma\T_{K_{\si{\infty}}}= \,\T_{K_{\si{\infty}}}/\,\T{\!}_{K_{\si{\infty}}}^{\;(\gamma-1)}$ n'est rien d'autre que le $\ell$-groupe des classes logarithmiques introduit dans \cite{J28} et calculé dans \cite{BJ} et \cite{DJ+}:\smallskip

\centerline{${}^\Gamma\T_{K_{\si{\infty}}}=\,\wCl_K$.}\smallskip

\item Et le sous-groupe ambige $\,\T{\!}_{K_{\si{\infty}}}^{\;\Gamma}$ est donné par un isomorphisme de Kuz'min (cf. \cite{Kuz}, Prop. 7.5 ou e.g. \cite{J55}, Th. 17, pour son interprétation logarithmique); il s'identifie au quotient\smallskip

\centerline{$\T{\!}_{K_{\si{\infty}}}^{\;\Gamma} \simeq \,\wE^\ph_K/\N_K$,}\smallskip

du 
groupe $\,\wE_K$ des unités logarithmiques de $K$ par le sous-groupe $\N_K
=\bigcap_{n\in\NN}N_{K_n/K}(\E'_{K_n})
$ des {\em normes universelles}, défini comme intersection des groupes de normes de $\ell$-unités,
mais qui est encore 
l'intersection $\,\wE_K^\nu=\bigcap_{n\in\NN}N_{K_n/K}(\wE_{K_n})$ des sous-groupes de normes d'unités logarithmiques  attachés aux divers étages de la tour cyclotomique (cf. Lem. \ref{NU} infra).
\end{itemize}\smallskip

En résumé, il vient ainsi:

\begin{Th}\label{TP}
Étant donné un corps de nombres totalement réel $K$ qui vérifie la conjecture de Leopoldt pour un premier donné $\ell$ (par exemple un corps abélien réel), la conjecture de Greenberg pour le corps $K$ et le premier $\ell$ affirme exactement que l'indice du sous-groupe des normes universelles  $\N_K=\,\wE_K^\nu=\bigcap_{n\in\NN}N_{K_n/K}(\wE_{K_n})$ dans le pro-$\ell$-groupe des unités logarithmiques $\,\wE_K$ du corps $K$ coïncide avec le cardinal du $\ell$-groupe des classes logarithmiques:\smallskip

\centerline{$(\,\wE_K^\ph : \N_K\,)\, =\,|\,\wCl_K^\ph|$.}
\end{Th}

Dans la formule obtenue, le $\ell$-groupe des classes logarithmiques $\,\wCl_K$ est effectivement calculable (cf. \cite{BJ,DJ+}) et un algorithme performant est désormais implanté dans {\sc pari} (cf. \cite{BJ}). Il en est de même du pro-$\ell$-groupe $\,\wE_K$, à cette réserve près que, les unités logarithmiques étant de nature $\ell$-adique, l'algorithme fournit un système minimal de générateurs de $\,\wE_K$ déterminés à puissance $\ell^{n_{\si{0}}}$-ième près, pour tout $n_{\si{0}}$ fixé à l'avance. En revanche, le sous-groupe des normes universelles, quoique bien défini dans $K$, reste difficilement calculable. La formule du Théorème fournit ainsi un critère (conjectural) d'arrêt. Donnons tout de suite un exemple (cf. Prop. \ref{Sco} infra):

\begin{Cor}
Soient $K$ un corps abélien réel de degré $d=[K:\QQ]$ étranger à $\ell$ et de groupe $\Delta$, puis $\Zl[\Delta]=\bigoplus_\varphi \Zl[\Delta]e_\varphi=\bigoplus_\varphi \ZZ_\varphi$ la décomposition semi-locale de son algèbre de Galois.\par
Le groupe des unités logarithmiques de $K$, regardé comme $\Zl[\Delta]$-module, est alors somme directe de ses composantes isotypiques $\,\wE_K^{\,e_\varphi}$, lesquelles sont $\ZZ_\varphi$-monogènes (à l'exception éventuelle de la composante unité $\{\pm1\}^{\Zl}\times \ell^{\Zl}$, qui est formée exclusivement de normes unverselles).\par
Sous la conjecture de Greenberg, la $\varphi$-composante isotypique $\N_K^{\,e_\varphi}$ du sous-groupe des normes universelles $\N_K$ est ainsi, pour chaque caractère $\varphi$, l'unique sous-module de $\,\E_K^{\,e_\varphi}$ d'indice $|\wCl_K^{\,e_\varphi}|$.
\end{Cor}

\newpage
\section{Interprétation en termes de capitulation logarithmique}

Considérons, comme dans le Théorème \ref{TP}, un corps totalement réel $K$ et un nombre premier $\ell$; notons $K_n$ le $n$-ième étage de la $\Zl$-tour cyclotomique $K_{\si{\infty}}/K$ et $\Gamma_n$ le groupe cyclique $\Gal(K_n/K)$.\smallskip


Le quotient $\,\wE^\ph_K/N_{K_n/K}(\wE_{K_n})$ s'identifie au groupe de cohomologie $H^2(\Gamma_n,\,\wE_K^\ph)$ relatif à l'action de $\Gamma_n=\Gal(K_n/K)$ sur le $\ell$-groupe des unités logarithmiques de $K_n$. Or, sous la conjecture de Gross-Kuz'min dans $K_n$, donc en particulier ici sous la conjecture de Leopoldt, le quotient de Herbrand $q(\Gamma_n,\,\wE_K^\ph)=|H^2(\Gamma_n,\,\wE_K^\ph)|/|H^1(\Gamma_n,\,\wE_K^\ph)|$ est égal à 1 (cf. \cite{J28}, Th. 3.6); de sorte que l'on a l'égalité:

\centerline{$(\,\wE^\ph_K:\,N_{K_n/K}(\wE_{K_n}))= |H^2(\Gamma_n,\,\wE_K^\ph)| = |H^1(\Gamma_n,\,\wE_K^\ph)|$.}\smallskip

D'autre part, l'extension $K_n/K$ étant logarithmiquement non ramifiée, le groupe $H^1(\Gamma_n,\,\wE_K^\ph)$ mesure la capitulation logarithmique dans $K_n/K$ (cf. \cite{J28}, Th. 4.5):\smallskip

\centerline{$H^1(\Gamma_n,\,\wE_K^\ph) \simeq \wCap_{K_n/K}=\Ker(\,\wCl_K\to\,\wCl_{K_n})$.}\smallskip

\noindent L'ordre de celle-ci étant borné par le cardinal du groupe des classes logarithmiques $\,\wCl_K$, la suite décroissante des groupes $N_{K_n/K}(\wE_{K_n}))$ stationne donc à partir d'un certain rang, disons, $n_{\si{0}}$:

\centerline{$\N_K=\wE_K^\nu=\bigcap_{n\in\NN}N_{K_n/K}(\wE_{K_n}^\ph)=N_{K_n/K}(\wE_{K_n}^\ph)$, pour tout $n\ge n_{\si{0}}$.}\smallskip

\noindent Et le résultat de capitulation donné dans \cite{J60, Ng} peut ainsi être précisé comme suit:

\begin{Th}
Sous la conjecture de Gross-Kuz'min dans $K_{\si{\infty}}$, l'indice des normes universelles $(\,\wE^\ph_K:\,\N_K^\ph\,)$ mesure la capitulation logarithmique dans $K_n/K$ pour tout $n$ assez grand:\smallskip

\centerline{ $(\,\wE^\ph_K:\N_K\,) = |\,\wCap_{K_{\si{\infty}}/K}|=  |\,\wCap_{K_n/K}|$, pour $n\gg 0$.}\smallskip

\noindent En particulier on a: $(\,\wE^\ph_K:\N_K\,) \le |\,\wCl_K|$ et l'égalité si et seulement si le corps $K$ vérifie la conjecture de Greenberg pour le premier $\ell$.
\end{Th}

Sous des hypothèses de semi-simplicité, il est facile, en outre, de raffiner l'inégalité obtenue:
supposons que le corps $K$ soit abélien sur un sous-corps $F$ avec $\ell\nmid [K:F]$ et notons $\Delta=\Gal(K/F)$ son groupe de Galois. Dans ce cas, l'algèbre de Galois $\Zl[\Delta]$ est un anneau semi local, produit direct d'extensions non-ramifiées $\ZZ_\varphi$ de $\Zl$ indexées par les caractères $\ell$-adiques irréductibles de $\Delta$; ce qui permet d'écrire tout $\Zl[\Delta]$-module $M$ comme somme directe de ses composantes isotypiques:\smallskip

 \centerline{$M=\oplus_\varphi M^{e_\varphi}$, avec $e_\varphi=\frac{1}{|\Delta|}\sum_{\sigma\in\Delta}\varphi(\sigma)\sigma^{\si{-1}}$.}

Il vient ainsi:

\begin{Prop}\label{Sco}
Sous les mêmes hypothèses et lorsque en outre le corps totalement réel $K$ est abélien sur un sous-corps $F$ de degré relatif $[K:F]$ étranger à $\ell$, l'égalité précédente vaut pour chaque composante isotypique de l'algèbre de Galois $\Zl[\Delta]$ associée au groupe $\Delta=\Gal(K/F)$:\smallskip

\centerline{ $(\,\wE^{\,e_\varphi}_K:\N_K^{\,e_\varphi}\,) = |\,\wCap_{K_n/K}^{\,e_\varphi}|$, pour $n\gg 0$.}\smallskip

Ce résultat vaut, en particulier, si $K$ est un corps abélien réel de degré $[K:\QQ]$ pour tout premier $\ell$ qui ne divise pas $[K:\QQ]$, auquel cas les $\varphi$-composantes du groupe des unités logarithmiques $\,\wE_K$ sont respectivement libres de dimension 1 sur l'anneau $Z_\varphi=\Zl[\Delta]e_\varphi$.
\end{Prop}

\Preuve Dans le cas abélien, en effet, le corps (totalement) réel $K$ vérifie la conjecture de Leopoldt pour tout premier $\ell$ et, par suite, celle de Gross-Kuz'min; laquelle assure que le caractère de $\,\wE_K$ est le caractère régulier (cf. \cite{J28}, Th. 3.6). Sous l'hypothèse $\ell\nmid[K:\QQ]$, c'est donc le produit de son sous groupe de torsion $\mu_K$ (qui est trivial pour $\ell\ne 2$) et d'un $\Zl[\Delta]$-module libre de dimension 1.

\Exemple Pour $\ell$ impair et $K$ quadratique réel, le groupe  $\,\wE_K$ des unités logarithmiques de $K$ est la somme directe de sa composante unité 
$\ell^\Zl$, qui est formée de normes universelles, et de sa composante d'augmentation, qui est un $\Zl$-module libre de rang 1. En particulier, le quotient $\,\wE^\ph_K/N_{K_n/K}(\,\wE_{K_n})$ est cyclique d'ordre inférieur ou égal au degré $\ell^n$ de l'extension $K_n/K$. Il suit de là que la capitulation logarithmique dans $K_n/K$ est d'ordre au plus $\ell^n$. Si le groupe des classes logarithmiques $\,\wCl_K$ est d'ordre $\ell^m$, il faut donc monter au moins jusqu'au $m$-ième étage de la tour cyclotomique pour le faire capituler, quel que soit l'exposant de $\,\wCl_K$.

\section{Comparaison avec la formule d'indice des unités circulaires}

La formule conjecturale donnée par le Théorème \ref{TP} peut être mise en parallèle avec l'identité\smallskip

\centerline{$(\,\E_K^\ph : \,\E_K^\circ\,)\, =\,|\,\Cl_K^\ph|$.}\smallskip

\noindent reliant l'indice du sous-groupe circulaire $\,\E_K^\circ$ dans le $\ell$-adifié $\,\E_K=\Zl\otimes_\ZZ E_K$ du groupe des unités au cardinal $|\,\Cl_K^\ph\,|$ du $\ell$-groupe des classes d'idéaux, obtenue par voie analytique (cf. \cite{Wa}, Th. 8.2).\smallskip

Donnons un exemple très simple:

\begin{Prop}
Soit $\ell$ un nombre premier impair, $n_{\si{0}}\ge 1$, puis $K=\QQ[\zeta+\zeta^{\si{-1}}]$ le sous-corps réel du corps cyclotomique $\QQ[\zeta]$ engendré par les racines $\ell^{n_{\si{0}}}$-ièmes de l'unité et $G=\Gal(K/\QQ)$ son groupe de Galois.\par

Le corps $K$ vérifie la conjecture de Greenberg pour $\ell$ si et seulement si le groupe $\N_K$ des normes universelles de $K$ coïncide avec le groupe circulaire $\,\C_K^\circ$, multiplicativement engendré sur l'algèbre de groupe $\Zl[G]$ par l'élément $\eta^\ph_K=(1-\zeta)(1-\zeta^{\si{-1}})$ construit sur une racine primitive $\ell^{n_{\si{0}}}$-ième de l'unité $\zeta$; autrement dit, si l'on a: $N_{K_n/K}(\E'_{K_n})=\,\C_K^\circ$ pour $n$ assez grand.
\end{Prop}

\Preuve
Rappelons que le $\ell$-groupe circulaire $\,\C^{\,\circ}_{K'}$ du corps cyclotomique $K'=\QQ[\zeta]$ est le $\Zl$-module multiplicatif engendré par la racine primitive $\zeta$, l'élément $\eta=1-\zeta$ et ses conjugués. Son sous-groupe réel $\,\C_K^\circ$ est l'image par la norme $N_{K'/K}=1+\tau$ ou, si l'on préfère ($\ell$ étant impair), par l'idempotent $e_+=\frac{1}{2}(1+\tau)$ associé à la conjugaison complexe $\tau$. C'est donc le $\Zl[G]$-module engendré par $\eta_K^\ph=\eta^{(1+\tau)}$, qui est contenu dans le $\ell$-adifié $\,\E'_K=\Zl\otimes_\ZZ E'_K$ du groupe des $\ell$-unités.\smallskip

Observons d'abord que toute $\ell$-unité $\varepsilon'\in\E'_K$ s'écrit de façon unique $\varepsilon'=\eta_K^\alpha\,\varepsilon$ avec $\alpha\in\Zl$ et $\varepsilon\in\E_K$; de sorte que nous avons la décomposition directe:\smallskip

\centerline{$\E'_K = \,\eta_K^\Zl \oplus \,\E_K$.}\smallskip

Regardons maintenant le sous-groupe circulaire $\,\C^\circ_K$. Il contient évidemment $\eta_K^\Zl$. Et son intersection avec le $\Zl[G]$-module des unités $\,\E_K$ est, par convention, le $\ell$-groupe des unités circulaires $\,\E_K^\circ$. Ainsi, de la décomposition

\centerline{$\C^\circ_K = \,\eta_K^\Zl \oplus \,\E^\circ_K$,}

\noindent nous concluons:

\centerline{$\E'_K/\,\C_K^\circ\simeq\,\E_K^\ph/\,\E_K^\circ$, \quad puis \quad $(\,\E'_K:\,\C_K^\circ\,)=(\,\E_K^\ph : \,\E_K^\circ\,) =|\,\Cl_K^\ph|$.}\smallskip

Cela étant, comme l'unique place de $K$ au-dessus de $\ell$ est principale, nous avons $\,\Cl_K^{\si{[\ell]}}=1$, i.e. $\,\Cl_K^\ph=\,\Cl_K'$. Portons enfin notre attention sur les $\ell$-groupes logarithmiques: rappelons qu'ils sont définis en remplaçant les valuations habituelles $v_\p$ attachées aux places finies de $K$ par leurs analogues logarithmiques $\wi v_\p$, lesquelles coïncident avec les précédentes pour $\p\nmid\ell$ et vérifient en outre une formule du produit (cf. \cite{J28}). Dans le cas qui nous occupe, puisque $K$ contient une unique place au-dessus de $\ell$ et que sa $\Zl$-extension cyclotomique $K_{\si{\infty}}/K$ est totalement ramifiée, nous avons directement: $\,\wCl_K^\ph=\,\Cl_K'$; ce qui nous donne: $(\,\E'_K:\,\C_K^\circ\,)=|\,\Cl'_K|=|\,\wCl_K^\ph|$.\smallskip

Enfin, toujours parce qu'il n'existe qu'une seule place au-dessus de $\ell$, le $\ell$-groupe $\,\wE_K$ des unités logarithmiques coïncide avec $\,\E'_K$. En fin de compte, l'égalité précédente s'écrit aussi bien:\smallskip

\centerline{ $(\,\wE_K:\,\C^\circ_K\,) =|\,\wCl_K^\ph|$;}\smallskip

\noindent tandis que la formulation de la conjecture de Greenberg donnée par le Théorème s'écrit, elle:\smallskip

\centerline{ $(\,\wE_K:\N_K\,) =|\,\wCl_K^\ph|$.}\smallskip

\noindent Et, comme on a l'inclusion banale $\,\C_K^\circ \subset \N_K$ du fait des identités normiques satisfaites par les unités circulaires $\eta^\ph_{\ell^{\si{n}}}=1-\zeta^\ph_{\ell^{\si{n}}}$, on voit que la conjecture de Greenberg pour $K$ et $\ell$ postule tout simplement l'égalité:

\centerline{$\,\C_K^\circ = \N_K$,}\smallskip

\noindent qui affirme que l'intersection $\bigcap_{n\in\NN}N_{K_n/K}(\E'_{K_n})$ se réduit à $\,\C_K^\circ$; autrement dit, que l'on a:\medskip

\centerline{$N_{K_n/K}(\E'_{K_n})=\,\C_K^\circ$, pour $n$ assez grand,}

\noindent comme annoncé.

\section{Étude du cas complètement décomposé}

On peut se demander si la coïncidence entre normes universelles et éléments circulaires donnée par la proposition précédente dans le cas particulier du sous-corps réel du corps cyclotomique $\QQ[\zeta^\ph_{\ell^n}]$ est générale. Il n'en est rien, comme le montre le cas complètement décomposé.\smallskip

Pour voir cela, partons d'un corps abélien réel $K$ dans lequel $\ell$ se décompose complètement; notons $G=G_K$ son groupe de Galois et $f=f_K$ son conducteur (qui n'est donc pas divisible par $\ell$). Puis, pour chaque diviseur $d>1$ de $f$, faisons choix d'une racine $d$-ième primitive de l'unité $\zeta_d^\ph$ et posons $\eta_d^\ph=1-\zeta_d^\ph$.\smallskip

Le pro-$\ell$-groupe des éléments circulaires (de Sinnott) $\,\C^\circ_K$ est alors le $\Zl[G]$-module multiplicatif engendré par le $\ell$-groupe $\mu_K^\ph$ des racines de l'unité contenues dans $K$ et les normes $N_{K[\zeta_d]/K}(\eta_d^\ph)$ des $\eta_d^\ph$. Il est bien connu que les $\eta_d^\ph$ sont des unités, pour $d$ non primaire; des $p$-unités, sinon. Dans tous les cas, ce sont donc des unités en $\ell$, de sorte que l'intersection de 
$\,\C^\circ_K$ avec le $\ell$-groupe des $\ell$-unités $\,\E'_K$ est contenu dans $\,\E_K^\ph$: c'est le $\ell$-groupe $\,\E^\circ_K$ des unités circulaires de Sinnott (cf. \cite{Si}).\smallskip

Or, comme observé par Greither (cf. \cite{Grt}, p. 218, l. 11--13), on a le résultat suivant:

\begin{Prop}
Pour tout corps abélien réel $K$ dans lequel le nombre premier $\ell$ est complètement décomposé l'intersection $\,\wE_K^\circ$ du $\ell$-groupe circulaire $\,\C^\circ_K$ avec le $\ell$-groupe de unités logarithmique $\,\wE_K^\ph$ se réduit au $\ell$-groupe $\mu_K^\ph$ des racines de l'unité contenues dans $K$:\smallskip

\centerline{ $\,\C_K^\circ \cap \,\wE_K^\ph = \,\E_K \cap \,\wE_K^\ph = \mu_K^\ph$.}
\end{Prop}

Posant $\,\wE_K^\circ=\C_K^\circ \cap \,\E_K$, on obtient donc, dans ce contexte:\smallskip

\centerline{$\wE_K^\circ=\mu_K^\ph$; \qquad mais \qquad $\N_K^\ph \simeq \mu_K^\ph \times \Zl^{[K:\QQ]}$,}\smallskip

\noindent puisque le groupe des normes universelles $\N_K$ est le produit de $\mu_K$ et d'un $\Zl$-module libre de rang $[K:\QQ]$, indépendamment de toute conjecture (cf. e.g. \cite{Grt}).\medskip

De fait, en l'absence même de toute hypothèse d'abélianité, l'indépendance (aux racines de l'unité près) entre groupes d'unités et groupes d'unités logarithmiques dans le cas complètement décomposé se lit déjà au niveau semi-local:

\begin{Th}
Soit $K$ un corps de nombres arbitraire dans lequel le nombre premier $\ell$ se décompose complètement; puis, pour chaque place $\l$ de $K$ au-dessus de $\ell$, soit $\ell_\l$ l'image canonique de $\ell$ dans le $\ell$-adifié $\R_{K_\l} =\varprojlim K_\l^\times/K_\l^{\times \ell^{\si{n}}}$ du groupe multiplicatif du complété $K_\l$.\par
Avec ces notations, le $\ell$-adifié $\R_{K_\l}$ s'écrit comme produit direct du sous-groupe des unités $\U_{K_\l}=\{\pm1\}^\Zl(1+\ell_\l)^\Zl$ et du sous-groupe des unités logarithmiques $\, \wU_{K_\l}=\ell_\l^\Zl$:\smallskip

\centerline{$\R_{K_\l}=\U_{K_\l}\,\wU_{K_\l}=\{\pm1\}^\Zl(1+\ell_\l)^\Zl \ell_\l^\Zl$.}
\end{Th}

\begin{Cor}
Lorsque $K$ est totalement réel, et sous la conjecture de Leopoldt pour $\ell$, l'application de semi-localisation $s_\ell$ identifie ainsi le produit $\,\E'_K\,(1+\ell)^\Zl$ du $\ell$-adifié $\,\E'_K$ du groupe des $\ell$-unités par le $\Zl$-module libre engendré par $1+\ell$ à un sous-module d'indice fini du produit:\smallskip

\centerline{ $\R_{K_\ell}=\prod_{\l\mid\ell}\R_{K_\l}=\prod_{\l\mid\ell}\{\pm1\}^\Zl(1+\ell_\l)^\Zl \ell_\l^\Zl$.}
\noindent Il suit:


\centerline{$\,\E_K\,(1+\ell)^\Zl \cap \,\wE_K = 1$, pour $\ell\ne 2$;\qquad $\,\E_K\,(1+\ell)^\Zl \cap \,\wE_K = \{\pm1\}$, pour $\ell=2$.}
\end{Cor}

\Preuve Dès lors que la place $\ell$ est complètement décomposée dans $K/\QQ$, chacun des complétés $K_\l$ au-dessus de $\ell$ s'identifie à $\Ql$. La décomposition des $\ell$-adifiés $\R_{K_\l}$ résulte donc directement de celle de $\QQ_\ell^\times$:

\centerline{$\QQ_\ell^\times =(1+\ell)^\Zl\ell^\ZZ$, pour $\ell$ impair;\qquad $\QQ_\ell^\times = \{\pm1\} (1+4)^\Zl 2^\ZZ$, pour $\ell=2$.}\smallskip

Si, de plus, $K$ vérifie la conjecture de Leopoldt pour le premier $\ell$ (cf. e.g. \cite{J31}, \S2.3), l'application de semi-localisation $s_\ell$ envoie injectivement le tensorisé $\,\E'_K=\Zl\otimes E'_K$ du groupe des $\ell$-unités de $K$ dans le produit $\R_{K_\ell}=\prod_{\l\mid\ell}\R_{K_\l}$ et se prolonge donc, du fait de l'égalité des rangs, en un pseudo-isomorphisme injectif du produit direct $\,\E'_K \, (1+\ell)^\Zl$ dans $\R_{K_\ell}$. D'où le résultat annoncé.

\section{Non-trivialité du quotient normique $\,\wE_K/\N_K$}

Revenons au cas général: considérons un corps totalement réel $K$ et intéressons-nous à la trivialité du quotient $\,\wE_K/\N_K$ sous la conjecture de Gross-Kuz'min. Notons $\,\T_{K_{\si{\infty}}} = \varprojlim \,\wCl_{K_n}$ le module de Kuz'min-Tate et écrivons $\,\F_{K_{\si{\infty}}}$ son sous-module fini. Il est bien connu que $\,\F_{K_{\si{\infty}}}$ est la limite projective des sous-groupes de capitulation $\wt{C\!ap}_{K_{\si{\infty}}/K_n}=\Ker\,(\wCl_{K_n}\to\wCl_{K_{\si{\infty}}})$ (cf. \cite{J24, J55,MN}:\smallskip

\centerline{$\,\F_{K_{\si{\infty}}} = \varprojlim \,\wt{C\!ap}_{K_{\si{\infty}}/K_n}$.}\smallskip

De l'inclusion $\T_{K_{\si{\infty}}}^{\,\Gamma} \subset\,\F_{K_{\si{\infty}}}$, donnée par la conjecture de Gross-Kuz'min, on tire directement:\smallskip

\centerline{$\T_{K_{\si{\infty}}}^{\,\Gamma}=1 \Rightarrow \,\F_{K_{\si{\infty}}}^{\,\Gamma}=1  \Rightarrow \,\F_{K_{\si{\infty}}}^\ph=1 \, {\Rightarrow} \,\T_{K_{\si{\infty}}}^{\,\Gamma}=1$.}\smallskip

\noindent de sorte que la trivialité de $\T{\!}_{K_{\si{\infty}}}^{\;\Gamma} \simeq \,\wE^\ph_K/\N_K$, équivaut à celle de  $\,\F_K$, ce qui revient à postuler que les morphismes d'extension $\,\wCl_{K_n}\to\wCl_{K_{\si{\infty}}}$ sont ultimement injectifs. Plus précisément:

\begin{Prop}
Sous la conjecture de Gross-Kuz'min dans $K_{\si{\infty}}$, on a $\,\wE_K=\N_K$, i.e. $\,\F_{K_{\si{\infty}}}=1$  si et seulement si le groupe des classes logarithmiques $\,\wCl_K$ s'injecte dans $\,\wCl_{K_{\si{\infty}}}$, i.e.  $\wt{C\!ap}_{K_{\si{\infty}}/K}=1$.
\end{Prop}

\Preuve D'après la suite exacte des classes logarithmiques ambiges (cf. \cite{J28}, Th. 4.5), il vient:\smallskip

\centerline{$\wE_K = \N_K  \underset{n\gg 0}{\Leftrightarrow} \,\wE_K = N_{K_n/K}(\wE_{K_n})  \underset{n\gg 0}{\Leftrightarrow} H^1(\Gamma_n,\wE_{K_n})=1  \underset{n\gg 0}{\Leftrightarrow} \,\wt{C\!ap}_{K_n/K}=1  \underset{n\gg 0}{\Leftrightarrow} \,\wt{C\!ap}_{K_{\si{\infty}}/K}=1$.}\smallskip

\noindent Ainsi, il suffit que  $\,\wCl_K$ s'injecte dans $\,\wCl_{K_{\si{\infty}}}$, pour qu'il en soit de même de tous les $\,\wCl_{K_n}$.\smallskip

\Remarque La condition $\,\wE_K/\N_K=1$ peut être regardée comme un critère suffisant de la conjecture de Gross-Kuz'min, a priori moins exigeant que la condition $\,\wCl_K=1$. Pour $K$ totalement réel cependant, la conjecture de Greenberg postule qu'elles sont équivalentes, en vertu du Théorème \ref{TP}.\medskip

Il peut être intéressant de donner un critère suffisant effectif de non-trivialité de $\,\wE_K / \N_K$.\smallskip

Partons d'un corps totalement réel $K$ , supposons pour simplifier $\ell\ne 2$ et notons $L=F[\zeta+\bar\zeta]$ le sous-corps réel du corps $N=F[\zeta]$ engendré sur $F$ par une racine $\ell$-ième primitive de l'unité $\zeta$. Écrivons enfin $|\mu^\ph_L|=\ell^m$ l'ordre du $\ell$-groupe des racines de l'unité dans $N$; notons $N_{\si{\infty}}= \bigcup N_n$ la $\Zl$-tour cyclotomique de $N$ (avec, donc, $N_n=K[\zeta^\ph_{\ell^{\si{m+n}}}$); et introduisons les radicaux universels ${\mathfrak R}_L=(\Ql/\Zl)\otimes_\ZZ L_\ph^\times \subset{\mathfrak R}_{N_{\si{\infty}}}=(\Ql/\Zl)\otimes_\ZZ N_{\si{\infty}}^\times$ (cf. \cite{J23}); puis considérons:\smallskip
\begin{itemize}
\item le sous-radical normique  ${}^\ph_{\ell^{\si{m}}}\!{\mathfrak N}_L=\{\ell^{\si{-m}}\otimes x\in {\mathfrak R}_L\,|\, x \in \N_L\}$, image des normes universelles;\smallskip 

\item le radical logarithmique ${}^\ph_{\ell^{\si{m}}}\!\we_L=\{\ell^{\si{-m}}\otimes x\in {\mathfrak R}_L\,|\, x \in \wE_L\} \subset (\Ql/\Zl)\otimes_\Zl\wE_L$;\smallskip

\item le radical initial ${}^\ph_{\ell^{\si{m}}}{\mathfrak Z}_L=\{\ell^{\si{-m}}\otimes x\in {\mathfrak R}_L\,|\, N_{\si{\infty}}[\sqrt[\ell^{\si{m}}]{x}]\subset Z\} =\Rad(Z/N_{\si{\infty}})$ attaché au compositum $Z$ des $\Zl$-extensions de $N$;\smallskip

\item et le noyau universel de Tate  ${}^\ph_{\ell^{\si{m}}}{\mathfrak U}_L=\{\ell^{\si{-m}}\otimes x\in {\mathfrak R}_L\,|\, \{\zeta_{\ell_\ph^{\si{m}}}^\ph,x\}=1 \;{\rm dans }\; K_2(N)\}$.\smallskip
\end{itemize}

Il est bien connu que  ${}^\ph_{\ell^{\si{m}}}{\mathfrak U}_L$ est un $(\ZZ/\ell^m\ZZ)$-module libre de dimension $[L:\QQ]$; qu'il en est de même de  ${}^\ph_{\ell^{\si{m}}}\!\we_L$ sous la conjecture de Gross-Kuz'min dans $N$, et de ${}^\ph_{\ell^{\si{m}}}{\mathfrak Z}_L$ sous celle de Leopoldt; qu'enfin ${}^\ph_{\ell^{\si{m}}}\!{\mathfrak N}_L$ est contenu dans chacun des trois précédents (cf. \cite{J24,J55,J57,MN,Seo1,Seo2,Seo3,Tat}. En particulier, il ne peut coïncider avec ${}^\ph_{\ell^{\si{m}}}\!\we_L$ que s'il coïncide avec les trois. Descendant alors par la norme (ou par l'idempotent 
$\frac{1}{[L:K]}N_{L/K})$ de l'algèbre de Galois $\Zl[\Gal(L/K)]$, nous en déduisons:

\begin{Th}[Critère de non-trivialité]
Sous les conjectures de Leopoldt et de Gross-Kuz'min dans $N$ (par exemple lorsque le corps $K$ est abélien sur $\QQ$), le quotient normique $\,\wE_K/\N_K$ ne peut être trivial que s'il y a coïncidence entre le radical logarithmique  ${}^\ph_{\ell^{\si{m}}}\!\we_K$, le radical initial attaché aux $\Zl$-extensions ${}^\ph_{\ell^{\si{m}}}{\mathfrak Z}_K$ et le noyau universel de Tate ${}^\ph_{\ell^{\si{m}}}{\mathfrak U}_K$.
\end{Th}

\begin{Sco}
Lorsque $K$ est un corps abélien réel de degré $d=[K:\QQ]$ étranger à $\ell$ et de groupe de Galois  $\Delta=\Gal(K/\QQ)$, le résultat vaut pour chaque composante isotypique de l'algèbre semi-locale $\Zl[\Delta]$; en d'autres termes, pour chaque caractère $\ell$-adique irréductible $\varphi$ de $\Delta$ on a l'implication:\smallskip

\centerline{$\F_K^{\,e_\varphi}=1 \quad \Leftrightarrow \quad \N^{\,e_\varphi}_K=\wE^{\,e_\varphi}_K\quad \Rightarrow \quad {}^\ph_\ell\we^{\,e_\varphi}_K={}^\ph_\ell{\mathfrak Z}^{\,e_\varphi}_K={}^\ph_\ell{\mathfrak U}^{\,e_\varphi}_K$.}
\end{Sco}

\newpage

\section*{Appendice: description du groupe des normes universelles}
\addcontentsline{toc}{section}{Appendice}

L'interprétation des normes universelles comme normes logarithmiques s'obtient comme suit:

\begin{Lem}\label{NU}
Soient $\ell$ un nombre premier, $K$ un corps de nombres et $K_\%=\bigcup_{n\in\NN}K_n$ sa $\Zl$-extension cyclotomique. Pour tout ensemble fini $S$ de places finies de $K$ contenant celles au-dessus de $\ell$, l'intersection $\N_K=\bigcap_{n\in\NN}N_{K_n/K}(\E^S_{K_n})$ des groupes normiques attachés aux $\ell$-adifiés $\,\E^S_{K_n}=\Zl\otimes E^S_{K_n}$ des groupes de $S$-unités est indépendante de $S$ et coïncide avec l'intersection $\,\wE^\nu_K=\bigcap_{n\in\NN}N_{K_n/K}(\wE_{K_n})$ des sous-groupes normiques d'unités logarithmiques.

\end{Lem}

\Preuve Partons d'un  $\eta^\ph_{\si{0}}$ de $\N_K$ et écrivons-le $N_{K_n/K}(\varepsilon_n)$ avec $\varepsilon_n\in\E^S_{K_n}$, pour tout $n\in\NN$. 
Observons que $\eta^\ph_{\si{0}}$ est une unité logarithmique (puisque norme à chaque étage fini $K_n/K$ de la $\Zl$-tour $K_\%/K$) et que, le groupe $\,\E_{K_1}^S$ étant compact, la suite $\big(N_{K_n/K_1}(\varepsilon_n)\big)_{n\ge 1}$ possède une valeur d'adhérence $\eta^\ph_{\si{1}}$.
De l'égalité $\eta^\ph_{\si{1}}=\lim N_{K_{n_{\si{k}}}/K_{\si{1}}}(\varepsilon_{n_{\si{k}}})$, pour une suite extraite convenable $(\varepsilon_{n_k})^\ph_{k\in\NN}$, les groupes normiques $N_{K_n/K_{\si{1}}}(\E^S_{K_n})$ étant fermés, on tire: 

\centerline{$\eta^\ph_{\si{1}}\in N_{K_{n_{\si{k}}}/K_{\si{1}}}(\E^S_{K_{n_{\si{k}}}})$ pour tout $n_k$; et finalement $\eta^\ph_{\si{1}}\in N_{K_{n}/K_{\si{1}}}(\E^S_{K_n})$ pour tout $n\ge 1$.}\smallskip

\noindent D'où le résultat par itération du procédé, puisqu'on a par construction $\eta^\ph_{\si{0}}=N_{K_{\si{1}}/K_{\si{0}}}(\eta^\ph_{\si{1}})$.


\def\refname{\normalsize{\sc  Références}}

\medskip\noindent
{\small
\begin{tabular}{l}
Institut de Mathématiques de Bordeaux,
Université de {\sc Bordeaux} \& CNRS \\
351 cours de la libération,
F-33405 {\sc Talence} Cedex\\
courriel : Jean-Francois.Jaulent@math.u-bordeaux.fr \\
\url{https://www.math.u-bordeaux.fr/~jjaulent/}
\end{tabular}
}

 \end{document}